\documentclass{article}
\usepackage[russian]{babel}
\usepackage{latexsym,amssymb}

\usepackage[cp1251]{inputenc}
\usepackage{amsfonts,amsmath}
\usepackage{graphicx,graphics,hhline}
\usepackage{euscript}
\newcounter{theorem} 
\newcounter{lemma} 
\renewcommand{\thetheorem}{\arabic{theorem}}

\newcommand{\theor}{\par\refstepcounter{theorem}%
{\bf Теорема \thetheorem .}\,\,}

\begin{document}

\centerline{\large\bf А.~К.~Бахтин} \vspace*{4mm}

\centerline{\large\textbf{\textit{Частично--конформные отображения в
}}}\vspace*{2mm} \centerline{\large\textbf{\textit{многомерных
 комплексных пространствах}}}


\vspace*{10mm}


В данной работе результаты статьи  \cite{32}  распространяются на
бесконечномерный случай. В своих исследованиях мы придерживаемся
терминологической архитектуры комплексного анализа, разработаного в
\cite{b1}-\cite{b5}.

Пусть $\mathbb{N, R, C}$ -- соответственно множества натуральных,
вещественных и комплексных чисел. $R_{+}=[0, +\infty)$. Пусть
$\mathbb{\overline{C}}$ -- сферa Римана (расширенная комплексная
плоскость), $r(B,a)$ -- внутренний радиус области
$B\subset\overline{C}$, относительно точки $a\in B$ (см. напр.
\cite{b6}-\cite{10}). По аналогии с пространством $\mathbb{C}^{n}$
рассмотрим линейное векторное пространство $\mathbb{C}^{\infty}$, то
есть пространство упорядоченных, счетных последовательностей
комплексных чисел. Таким образом $\mathbb{C}^{\infty}$ есть
декартовое произведение счетного числа экземпляров комплексной
плоскости $\mathbb{C}$:
$\mathbb{C}^{\infty}=\mathbb{C}\times\mathbb{C}\times\ldots\times\mathbb{C}\times\ldots$.
Аналогично,
$\mathbb{R}^{\infty}=\mathbb{R}\times\mathbb{R}\times\ldots\times\mathbb{R}\times\ldots$,
$\mathbb{R}^{\infty}\subset\mathbb{C}^{\infty}$.

Сформулируем для случая пространства $\mathbb{C}^{\infty}$ некоторые
понятия работы \cite{32}.

{\bf 1. Алгебра $\mathbb{C}^{\infty}$.}

\textbf{\emph{Определение 1.}} Бинарную операцию действующую из
$\mathbb{C}^{\infty}\times\mathbb{C}^{\infty}$ в
$\mathbb{C}^{\infty}$ по правилу
\begin{equation}
\label{tri} \mathbb{Z}\cdot\mathbb{W}=
\{z_{k}{w}_{k}\}_{k=1}^{\infty},
\end{equation}
где $\mathbb{Z}=\{z_{k}\}_{k=1}^{\infty}\in \mathbb{C}^{\infty}$,
$\mathbb{W}=\{w_{k}\}_{k=1}^{\infty}\in \mathbb{C}^{\infty}$, будем
называть \emph{векторным умножением элементов
$\mathbb{C}^{\infty}$}. Данная операция превращает
$\mathbb{C}^{\infty}$ в коммутативную, ассоциативную алгебру
\cite{b7, b8} с единицей
$\textbf{\textsf{1}}=(1,1,...,1,\ldots)\in\mathbb{C}^{\infty}$.

Обратимыми, относительно так определенной операции умножения,
являются те и только те елементы
$\mathbb{Z}=\{z_{k}\}_{k=1}^{\infty}\in \mathbb{C}^{\infty}$ у
которых $z_{k}\neq0$ для всех $k=\overline{1,\infty}$.

Обратными для таких элементов $\mathbb{Z}\in \mathbb{C}^{\infty}$
являются элементы $\mathbb{Z}^{-1}=\{z_{k}^{-1}\}_{k=1}^{\infty}\in
\mathbb{C}^{\infty}$, так как
$\mathbb{Z}\cdot\mathbb{Z}^{-1}=\mathbb{Z}^{-1}\cdot\mathbb{Z}=1$.
Множество $\Theta$ всех элементов
$A=\{a_{k}\}_{k=1}^{\infty}\in\mathbb{C}^{\infty}$, у которых хотя
бы одна координата $a_{k}=0$, назовем множеством необратимых
элементов $A\in \mathbb{C}^{\infty}$. Множество $\Theta$ является
объединением максимальных идеалов алгебры $\mathbb{C}^{\infty}$
\cite{b10}.

По аналогии с конечномерным случаем (см. напр \cite{b7},
\cite{b8}), можно представить $\mathbb{C}^{\infty}$ как
прямую сумму счетного числа экземпляров алгебры комплексных чисел
$\mathbb{C}$. Структура векторного пространства
$\mathbb{C}^{\infty}$ полностью согласуется со структурой алгебры
$\mathbb{C}^{\infty}$. Обобщим ряд определений работы \cite{32} на
случай алгебры $\mathbb{C}^{\infty}$.

{\bf 2. Сопряжение.}

\textbf{\emph{Определение 2.}} Каждому элементу
$\mathbb{W}=\{w_{k}\}_{k=1}^{\infty}\in \mathbb{C}^{\infty}$
поставим в соответствие векторно -- сопряженный элемент
$\overline{\mathbb{W}}=\{\overline{w}_{k}\}_{k=1}^{\infty}\in
\mathbb{C}^{\infty}$, где $\overline{w}_{k}$ обозначает число
комплексно сопряженное $w_{k}$ в обычном смысле. Так определенное
соответствие задает автоморфизм $\mathbb{C}^{\infty}$, оставляющий
неподвижным подпространство $\mathbb{R}^{\infty}$.

{\bf 3. Модуль (векторный).} В алгебре $\mathbb{C}$ одним из
важнейших является понятие модуля комплексного числа. В работе
\cite{32} предложено векторное обобщение даного понятия,
распространим его на пространство $\mathbb{C}^{\infty}$. Пусть
$\mathbb{R}_{+}^{\infty}=R_{+}\times R_{+}\times \ldots \times
R_{+}\ldots$.

\textbf{\emph{Определение 3.}} Векторным модулем произвольного
элемента $\mathbb{Z}=\{z_{k}\}_{k=1}^{\infty}\in
\mathbb{C}^{\infty}$ будем называть вектор
$|\mathbb{Z}|:=\{|z_{k}|\}_{k=1}^{\infty}\in
\mathbb{R}_{+}^{\infty}$.

Важно, что для произвольного $\mathbb{Z}=\{z_{k}\}_{k=1}^{\infty}\in
\mathbb{C}^{\infty}$, справедливо равенство
\begin{equation}\label{triv}
\mathbb{Z}\cdot\overline{\mathbb{Z}}=
|\overline{\mathbb{Z}}|^{2}=|\mathbb{Z}|^{2}.
\end{equation}

{\bf 4. Векторная норма.}

\textbf{\emph{Определение 4.}} Вектор
$\mathbb{X}=\{x_{k}\}_{k=1}^{\infty}\in \mathbb{R}^{\infty}$ будем
называть неотрицательным (строго положительным) и писать
$\mathbb{X}\geq \mathbb{O}$ ($\mathbb{X}> \mathbb{O}$), если
$x_{k}\geq 0$  для всех $k=\overline{1, \infty}$ ($x_{k}>0$ хотя бы
для одного $k=\overline{1, \infty}$), $\mathbb{O}=(0, 0, \ldots ,0,
\ldots)$.

\textbf{\emph{Определение 5.}} Будем говорить, что вектор
$\mathbb{X}=\{x_{k}\}_{k=1}^{\infty}\in \mathbb{R}^{\infty}$ больше
либо равен (строго больше) вектора
$\mathbb{Y}=\{y_{k}\}_{k=1}^{\infty}\in \mathbb{R}^{\infty}$, если
$\mathbb{X}-\mathbb{Y}\geq\mathbb{O}$
($\mathbb{X}-\mathbb{Y}>\mathbb{O}$).

Данные определения для конечномерных пространств $\mathbb{C}^{n}$
были даны в работе \cite{32}. В многомерных пространствах ситуация
существенно отличается от случая вещественной прямой, например,
вектор $\mathbb{O}=\underbrace{(0, 0, \ldots 0,
\ldots)}_{\infty-\textrm{раз }}$ больше либо равен всех векторов,
координаты которых неположительны и меньше либо равен всех векторов
из $\mathbb{R}_{+}^{\infty}$. Остальные векторы
$\mathbb{R}^{\infty}$ у которых координаты разных знаков с вектором
$\mathbb{O}$ не сравнимы в смысле определений 4 и 5.

\textbf{\emph{Определение 6.}} Векторное пространство $\mathbb{Y}$
будем называть векторно нормированным, если каждому $y\in\mathbb{Y}$
сопоставлен неотрицательный вектор
$\|y\|\in\mathbb{R}_{+}^{\infty}$, удовлетворяющий условиям:

1) $\|y\|\geq\mathbb{O}$, причем
$\|y\|=\mathbb{O}\Longleftrightarrow y=0_{\mathbb{Y}}$,
($0_{\mathbb{Y}}$ -- нуль пространства $\mathbb{Y}$);

2) $\|\gamma y\|=|\gamma |\|y\|$, $\forall y\in\mathbb{Y}$,
$\forall\gamma\in\mathbb{C}$;

3) $\|y_{1}+y_{2}\|\leq\|y_{1}\|+\|y_{2}\|$, $\forall y_{1},
y_{2}\in\mathbb{Y}$. Аналогично можно ввести понятие векторной
метрики. Введенное определение 3 удовлетворяет определению 6. Таким
образом векторный модуль является векторной нормой в алгебре
$\mathbb{C}^{\infty}:\|\cdot\|=|\cdot|$. Тогда открытым единичным
шаром в алгебре $\mathbb{C}^{\infty}$ является единичный открытый
поликруг $\|z\|<\textbf{\textsf{1}}$, $(\textbf{\textsf{1}}=(1, 1,
\ldots, 1, \ldots)$, а единичной сферой --
$\mathbb{T}^{\infty}=\{\mathbb{Z}\in\mathbb{C}^{\infty}:\|Z\|=\textbf{\textsf{1}}\}$.
Важно, что

а) $|Z_{1}\cdot Z_{2}|=\|Z_{1}\cdot
Z_{2}\|=\|Z_{1}\|\|Z_{2}\|=|Z_{1}||Z_{2}|$, $\forall Z_{1},
Z_{2}\in\mathbb{C}^{\infty}$;

б) $|1|=\|1\|=\textbf{\textsf{1}}$, $(\textbf{\textsf{1}}=(1,
1,\ldots, 1, \ldots)).$

Заметим, что для евклидовой нормы $\|\cdot\|_{E}$, определяемой
соотношением
$\|z\|_{E}=\sum\limits_{k=1}^{\infty}(|z_{k}|^{2})^{\frac{1}{2}}$,
справедливо равенство
$$\|1\|_{E}=\infty.$$

{\bf 5. Векторный аргумент $a\in\mathbb{C}^{\infty}$.} В дальнейшем
вектор (произвольный) пространства (алгебры) $\mathbb{C}^{\infty}$
будем называть  бесконечномерным комплексным числом, а алгебра
$\mathbb{C}^{\infty}$ будет называться алгеброй бесконечномерных
комплексных чисел.

\textbf{\emph{Определение 7.}} Векторным аргументом
бесконечномерного комплексного числа
$\mathbb{A}=\{a_{k}\}_{k=1}^{\infty}\in
\mathbb{C}^{\infty}\backslash\Theta$ является бесконечномерный
вещественный вектор, определяемый формулой $$ \arg \mathbb{A}=\{\arg
a_{k}\}_{k=1}^{\infty},$$ где $\arg a_{k}$ есть главное значение
аргумента, либо то которое вытекает из конкретного смысла задачи в
которой фигурирует бесконечномерное комплексное число
$\mathbb{A}\in\mathbb{C}^{\infty}$.

{\bf 6. Компактификация $\mathbb{C}^{\infty}$.}

В качестве компактификации
$\mathbb{C}^{\infty}=\mathbb{C}\times\mathbb{C}\times\ldots
\times\mathbb{C}\times\ldots$ возьмем пространство
$\overline{\mathbb{C}}^{\infty}=\overline{\mathbb{C}}\times\overline{\mathbb{C}}
\times\ldots \times\overline{\mathbb{C}}\times\ldots$, которое по
аналогии с конечномерным случаем (см. \cite{b1}-\cite{b5}) будем
называть бесконечномерным пространством теории функций. Бесконечными
точками $\overline{\mathbb{C}}^{\infty}$ являются те точки, у
которых хотя бы одна координата бесконечна. Множество всех
бесконечных точек имеет коразмерность единица.

Топологию в $\overline{\mathbb{C}}^{\infty}$ определяем как
покоординатную сходимость равномерную по номерам координат.

\textbf{7. Дифференцируемость.} Сначала обратимся к конечномерному
случаю. Рассмотрим область $\mathbb{D}\subset\mathbb{C}^{n}$ и
отображение $\mathbb{F}:\mathbb{D}\rightarrow \mathbb{C}^{m}$,
$\mathbb{F}=\{f_{k}(z_{1},\ldots,z_{n})\}_{k=1}^{m}.$ Пусть
$f_{k}=U_{k}(x_{1},\ldots,x_{n}, y_{1},\ldots,
y_{n})+iV_{k}(x_{1},\ldots,x_{n}, y_{1},\ldots, y_{n})$ --
вещественно непрерывно дифференцируемы всюду в области $\mathbb{D}$
при $k=\overline{1,m}$, $n,m \in \mathbb{N}$.

Матрицу Якоби отображения $\mathbb{F}$, рассматриваемого как
дифференцируемое отображение области
$\mathbb{D}\subset\mathbb{R}^{2n}$ в $\mathbb{R}^{2m}$ (матрица
$2m\times2n$) представим следующим образом

\begin{equation}\label{f3.6}\begin{pmatrix}
U_{x_{1}}^{(1)}&\ldots&U_{x_{n}}^{(1)}&|&U_{y_{1}}^{(1)}&\ldots&U_{y_{n}}^{(1)}\\
\hdotsfor{3}&|&\hdotsfor{3}\\
\vdots&\{\mathbb{U}_{\mathbb{X}}\}&\vdots&|&\vdots&\{\mathbb{U}_{\mathbb{Y}}\}&\vdots\\
\hdotsfor{3}&|&\hdotsfor{3}\\
U_{x_{1}}^{(m)}&\ldots& U_{x_{n}}^{(m)}&|&U_{y_{1}}^{(m)}&\ldots&U_{y_{n}}^{(m)}\\
---&---&---&|&---&---&---\\
V_{x_{1}}^{(1)}&\ldots&V_{x_{n}}^{(1)}&|&V_{y_{1}}^{(1)}&\ldots&V_{y_{n}}^{(1)}\\
\hdotsfor{3}&|&\hdotsfor{3}\\
\vdots&\{\mathbb{V}_{\mathbb{X}}\}&\vdots&|&\vdots&\{\mathbb{V}_{\mathbb{Y}}\}&\vdots\\
\hdotsfor{3}&|&\hdotsfor{3}\\
V_{x_{1}}^{(m)}&\ldots& V_{x_{n}}^{(m)}&|&V_{y_{1}}^{(m)}&\ldots&V_{y_{n}}^{(m)}\\
\end{pmatrix},\end{equation}

где $U_{x_{j}}^{(k)}=\frac{\partial}{\partial x_{j}}U_{k},$
$V_{x_{j}}^{(k)}=\frac{\partial}{\partial x_{j}}V_{k},$
$k=\overline{1,m}$, $j=\overline{1,n}$.

Штрихованные линии разбивают матрицу Якоби (\ref{f3.6}) на четыре
прямоугольные матрицы порядка $m\times n$, обозначенные
$\mathbb{U}_{\mathbb{X}}$, $\mathbb{U}_{\mathbb{Y}}$,
$\mathbb{V}_{\mathbb{X}}$, $\mathbb{V}_{\mathbb{Y}}$, где
$\mathbb{F}=Re \mathbb{F} +iIm\mathbb{F}=\mathbb{U}+i\mathbb{V},$
$\mathbb{Z}=Re\mathbb{Z}+iIm\mathbb{Z}=\mathbb{X}+i\mathbb{Y}.$

С учетом сказанного, матрицу (\ref{f3.6}) можно представить
следующим образом

\begin{equation}\label{f4.6}
\begin{pmatrix}
\mathbb{U}_{\mathbb{X}}& \mathbb{U}_{\mathbb{Y}}\\
\mathbb{V}_{\mathbb{X}}& \mathbb{V}_{\mathbb{Y}}\\
\end{pmatrix},\end{equation}

Тогда условия Коши-Римана для отображения $\mathbb{F}$ можно
записать в виде

\begin{equation}\label{f3.7}
\left\{
\begin{aligned}
\mathbb{U}_{\mathbb{X}}=\mathbb{V}_{\mathbb{Y}},\\
\mathbb{U}_{\mathbb{Y}}=-\mathbb{V}_{\mathbb{X}}.\\
\end{aligned}
\right.\end{equation}

С учетом (\ref{f3.7}) известное определение голоморфного отображения
(см. \cite{b1}-\cite{b5}) можно представить в следующем виде.

\textbf{\emph{Определение 8.}} Отображение
$\mathbb{F}:\mathbb{D}\rightarrow \mathbb{C}^{m}$ вещественно
непрерывно дифференцируемое в $\mathbb{D}$ (как отображение из
$\mathbb{R}^{2n}$ в $\mathbb{R}^{2m}$) и удовлетворяющее матричному
уравнению (\ref{f3.7}) всюду в $\mathbb{D}$ будем называть
голоморфным в области $\mathbb{D}$. При $n \in \mathbb{N}$ и $m=1$
получаем определение голоморфной функции в области
$\mathbb{D}\subset\mathbb{C}^{n}$. В случае $n=1$, $m \in
\mathbb{N}$ получаем определение голоморфной кривой.

Как известно \cite{b1}-\cite{b5} голоморфное отображение
$\mathbb{F}:\mathbb{D}\rightarrow \mathbb{C}^{m}$,
$\mathbb{D}\subset\mathbb{C}^{n}$ называется биголоморфным, если оно
имеет обратное отображение, голоморфное в области
$\mathbb{F}(\mathbb{D})$.

Теперь дадим формальное обобщение приведенных выше рассуждений на
бесконечномерный случай.

Пусть даны область $\mathbb{D}\subset\mathbb{C}^{\infty}$ и
отображение $\mathbb{F}:\mathbb{D}\rightarrow\mathbb{C}^{\infty}$,
где
$\mathbb{F}=\{f_{k}(\mathbb{Z})\}_{k=1}^{n}=\{f_{k}(\mathbb{X}+i\mathbb{Y})\}_{k=1}^{n}$,
$f_{k}(\mathbb{X}+i\mathbb{Y})=U_{k}(\mathbb{X},
\mathbb{Y})+iV_{k}(\mathbb{X},\mathbb{Y})=U_{k}(\{x_{p}\}_{p=1}^{\infty},\{y_{p}\}_{p=1}^{\infty})+iV_{k}(\{x_{p}\}_{p=1}^{\infty},\{y_{p}\}_{p=1}^{\infty})$.
$\mathbb{F}=\mathbb{U}+i\mathbb{V}$,
$\mathbb{U}=\mathbb{U}(\mathbb{X},\mathbb{Y})=\{U_{k}(\mathbb{X},
\mathbb{Y})\}_{k=1}^{\infty}$,
$\mathbb{V}=\mathbb{V}(\mathbb{X},\mathbb{\mathbb{Y}})=\{V_{k}(\mathbb{X},
\mathbb{Y})\}_{k=1}^{\infty}$,
$\mathbb{Z}=\mathbb{X}+i\mathbb{Y}=\{x_{k}\}_{k=1}^{\infty}+i\{y_{k}\}_{k=1}^{\infty}\in
\mathbb{D}$.

Пусть функции
$U_{k}(\{x_{p}\}_{p=1}^{\infty},\{y_{p}\}_{p=1}^{\infty})$,
$V_{k}(\{x_{p}\}_{p=1}^{\infty},\{y_{p}\}_{p=1}^{\infty})$ всюду в
$\mathbb{D}$ имеют непрерывные частные производные по всем
переменным $x_{p}$, $y_{p}$, $p=\overline{1,\infty}$. Тогда матрицу
Якоби представим в виде аналогичном (\ref{f4.6})

\begin{equation}
\begin{pmatrix}
\mathbb{U}_{\mathbb{X}}& \mathbb{U}_{\mathbb{Y}}\\
\mathbb{V}_{\mathbb{X}}& \mathbb{V}_{\mathbb{Y}}\\
\end{pmatrix},\end{equation}
где $\mathbb{U}_{\mathbb{X}}, \mathbb{U}_{\mathbb{Y}},
\mathbb{V}_{\mathbb{X}}, \mathbb{V}_{\mathbb{Y}}$ являются
бесконечными матрицами следующего вида
$\mathbb{U}_{\mathbb{X}}=\left[\{U_{x_{p}}^{(k)}\}_{k=1,
p=1}^{\infty}\right]$,
$\mathbb{U}_{\mathbb{Y}}=\left[\{U_{y_{p}}^{(k)}\}_{k=1,
p=1}^{\infty}\right]$,
$\mathbb{V}_{\mathbb{X}}=\left[\{V_{x_{p}}^{(k)}\}_{k=1,
p=1}^{\infty}\right]$,
$\mathbb{V}_{\mathbb{Y}}=\left[\{V_{y_{p}}^{(k)}\}_{k=1,
p=1}^{\infty}\right]$, $V_{x_{p}}^{(k)}=\frac{\partial}{\partial
x_{p}}V_{k}$, $V_{y_{p}}^{(k)}=\frac{\partial}{\partial
y_{p}}V_{k}$, $U_{x_{p}}^{(k)}=\frac{\partial}{\partial
x_{p}}U_{k}$, $U_{y_{p}}^{(k)}=\frac{\partial}{\partial
x_{p}}U_{k}$, $k,p=\overline{1,\infty}$.

В нашем случае, символ $[\cdot]$ обозначает бесконечную матрицу.

Тогда уравнения Коши-Римана примут следующий вид

\begin{equation}\label{f3.8}
\left\{
\begin{aligned}
\mathbb{U}_{\mathbb{X}}=\mathbb{V}_{\mathbb{Y}},\\
\mathbb{U}_{\mathbb{Y}}=-\mathbb{V}_{\mathbb{X}}.\\
\end{aligned}
\right.\end{equation}

Формально системы (\ref{f3.7}) и (\ref{f3.8}) совершенно одинаковы.

\textbf{\emph{Определение 9.}} Пусть $\mathbb{D}$ является
произвольной областью из пространства  $\mathbb{C}^{\infty}$.
Отображение $\mathbb{F}:\mathbb{D}\rightarrow \mathbb{C}^{\infty}$
вещественно непрерывно дифференцируемое в $\mathbb{D}$ и
удовлетворяющее матричному уравнению (\ref{f3.8}) всюду в
$\mathbb{D}$ будем называть голоморфным отображением области
$\mathbb{D}$.

Если $k\geq m$, $f_{k}(\mathbb{Z})\equiv 0$, то мы получаем
голоморфное отображение $\mathbb{F}:\mathbb{D}\rightarrow
\mathbb{C}^{m}$. Если рассмотрим сужение отображения $\mathbb{F}$ на
$\mathbb{C}^{n}$, то получим голоморфное отображение из
$\mathbb{C}^{n}$ в $\mathbb{C}^{\infty}$.

По аналогии с конечномерным случаем, будем считать, что голоморфное
отображение $\mathbb{F}:\mathbb{D}\rightarrow \mathbb{C}^{\infty}$,
$\mathbb{D}\subset\mathbb{C}^{\infty}$ является биголоморфным, если
$\mathbb{F}$ имеет обратное отображение, голоморфное в
$\mathbb{F}(\mathbb{D})$.

Для того чтоб перенести известные результаты о сходимости степенных
рядов, известные для комплексной плоскости $\mathbb{C}$, в
пространство $\mathbb{C}^{\infty}$, приведем соответствующий аналог
определения равномерной сходимости внутри единичного поликруга
некоторой последовательности функций отображений.

Пусть $\mathbb{U}_{r}^{\infty}=U_{r}\times U_{r}\times\ldots\times
U_{r}\times\ldots$, где $U_{r}=\{z:z\in \mathbb{C}, |z|<r\}$,
$\mathbb{U}_{1}^{\infty}:=\mathbb{U}^{\infty}$.
$\overline{\mathbb{U}}_{r}^{\infty}=\overline{U}_{r}\times
\overline{U}_{r}\times\ldots\times \overline{U}_{r}\times\ldots$, и
$\mathbb{F}_{p}:\mathbb{U}^{\infty}\rightarrow \mathbb{C}^{\infty}$
-- некоторая последовательность отображений.

\textbf{\emph{Определение 10.}} Будем говорить, что
последовательность $\mathbb{F}_{p}$, $p=\overline{1,\infty}$
равномерно внутри $\mathbb{U}^{\infty}$ сходится к некоторому
отображению $\mathbb{F}_{0}:\mathbb{U}^{\infty}\rightarrow
\mathbb{C}^{\infty}$, если для любого $\varepsilon>0$ и $0<r<1$
существует такой номер $n_{0}=n_{0}(\varepsilon,r)$, $n_{0}\in
\mathbb{N}$, что
$$\|\mathbb{F}_{p}(\mathbb{Z})-\mathbb{F}_{0}(\mathbb{Z})\|\leq\varepsilon\cdot\textbf{1}$$
для всех $\mathbb{Z}\in\overline{\mathbb{U}}_{r}^{\infty}$ и всех
$p>n_{0}$.

\textbf{\emph{Определение 11.}} Голоморфное отображение
$$\mathbb{F}:\mathbb{U}^{\infty}\rightarrow
\mathbb{C}^{\infty},\quad\mathbb{F}(\mathbb{Z})=\begin{pmatrix}f_{1}(z_{1})\\
f_{2}(z_{2})\\
\vdots\\
f_{n}(z_{n})\\
\vdots\\
\end{pmatrix},\quad f_{k}=\sum \limits_{p=0}^\infty a_{p}^{(k)}z_{k}^{p},$$ будем называть аналитической функцией
векторного аргумента, если равномерно внутри поликруга
$\mathbb{U}^{\infty}$ сходится ряд

$$\mathbb{F}(\mathbb{Z})=\sum \limits_{p=1}^\infty
\mathbb{A}_{p}\mathbb{Z}^{p},\quad{\mathbb{A}}_{p}=\{a_{p}^{(k)}\},\quad
\mathbb{Z}\in \mathbb{U}^{\infty},\quad p,k=\overline{1,\infty}.$$

Аналитические функции векторного аргумента подводят нас к понятию
\emph{частично -- конформных отображений} в $\mathbb{C}^{n}$ и
$\mathbb{C}^{\infty}$.

\textbf{\emph{Определение 12.}} Пусть $\delta>0$ -- некоторое
фиксированное
число. Тогда отображение $$\mathbb{F}(\mathbb{Z})=\begin{pmatrix}f_{1}(z_{1})\\
f_{2}(z_{2})\\
\vdots\\
f_{n}(z_{n})\\
\vdots\\
\end{pmatrix},\quad\mathbb{Z}\in\mathbb{U}^{\infty},$$ где каждое
$f_{k}(z_{k})$, $k=\overline{1,\infty}$, является однолистной
функцией в единичном круге, такой что
$\delta<|f'_{k}(0)|<\frac{1}{\delta}$, $k=\overline{1,\infty}$,
будем называть частично конформным отображением единичного
поликруга.

Поясним, что в определении 12 число $\delta$ зависит, вообще говоря,
от отображения $\mathbb{F}$, то есть $\delta=\delta(\mathbb{F})$.

Текже из определения 12 следует определение частично -- конформного
отображения в любом пространстве $\mathbb{C}^{n}$.

В целом, частично -- конформные отображения, вообще говоря, не
являются конформными отображениями, однако сужение таких отображений
на любую координатную плоскость есть конформным отображением.

{\bf 8. Представление бесконечномерного комплексного числа в
векторно -- декартовой формe.} Пусть
$\mathbb{Z}=\{z_{k}\}_{k=1}^{\infty}\in\mathbb{C}^{\infty}$. Тогда
$$\mathbb{Z}=\{z_{k}\}_{k=1}^{\infty}=\{Re z_{k}+i Im
z_{k}\}_{k=1}^{\infty}=\{Re z_{k}\}_{k=1}^{\infty}+\{i Im
z_{k}\}_{k=1}^{\infty}=$$
$$=\{Re z_{k}\}_{k=1}^{\infty}+i\{ Im z_{k}\}_{k=1}^{\infty}= Re\mathbb{Z}+i Im \mathbb{Z}=X+iY=$$
$$=\{x_{k}\}_{k=1}^{\infty}+i\{y_{k}\}_{k=1}^{\infty}\in\mathbb{R}^{\infty}+i\mathbb{R}^{\infty},$$
где $X=Re\mathbb{Z}=\{Re
z_{k}\}_{k=1}^{\infty}=\{x_{k}\}_{k=1}^{\infty},$
$Y=Im\mathbb{Z}=\{Im
z_{k}\}_{k=1}^{\infty}=\{y_{k}\}_{k=1}^{\infty}.$ То есть
$\mathbb{C}^{\infty}=\mathbb{R}^{\infty}+i\mathbb{R}^{\infty}$.

{\bf 9. Представление бесконечномерного комплексного числа в
векторно -- полярной форме.}

Используя вышеприведенные определения, получим цепочку равенств:

$$\mathbb{Z}=\{z_{k}\}_{k=1}^{\infty}=\{|z_{k}|e^{i\alpha_{k}}\}_{k=1}^{\infty}=\{|z_{k}|\}_{k=1}^{\infty}\{e^{i\alpha_{k}}\}_{k=1}^{\infty}$$
$$=|\mathbb{Z}|\left[\cos\arg\mathbb{Z}+i\sin\arg
\mathbb{Z}\right]=|\mathbb{Z}|e^{i\arg\mathbb{Z}},$$

где $$\cos\beta=\{\cos\beta_{k}\}_{k=1}^{\infty},\quad
\sin\beta=\{\sin\beta_{k}\}_{k=1}^{\infty},$$
$$\exp i\mathbb{\beta}=\{\exp i\beta_{k}\}_{k=1}^{\infty}, \quad \mathbb{\beta} =\{\beta_{k}\}_{k=1}^{\infty}\in \mathbb{R}^{\infty},\quad \mathbb{Z}=
\{z_{k}\}_{k=1}^{\infty}\in \mathbb{C}^{\infty}.$$

Аналогичным образом определяется отображение $\ln \mathbb{Z}$,\quad
$\mathbb{Z}=\{z_{k}\}_{k=1}^{\infty}\in
\mathbb{C}^{\infty}\setminus\Theta$
$$\ln \mathbb{Z}=\ln |\mathbb{Z}|+i\arg\mathbb{Z}=
\{\ln|z_{k}|+i\arg z_{k}\}_{k=1}^{\infty}.$$

Более того, для регулярной в областях $(B_{1}, B_{2},\ldots,B_{n},
\ldots)$, $B_{k}\in \mathbb{C}$, $k=\overline{1,\infty}$ функции
$F(z)$ комплексного переменного определим продолжение этой функции
до голоморфного отображения области $\mathbb{B}=B_{1}\times
B_{2}\times \ldots\times B_{n}\times\ldots$ по следующему правилу
$$\mathbb{F}(\mathbb{W})=\{F(W_{k})\}_{k=1}^{\infty},\quad
\mathbb{W}=\{w_{k}\}_{k=1}^{\infty}\in \mathbb{B}.$$

Приведем примеры частично -- конформных отображений, заданих
элементарными функциями.

\textbf{1)} Дробно -- линейная функция определяется соотношением
$$T=\frac{\mathbb{A}_{1}\mathbb{Z}+\mathbb{A}_{2}}{\mathbb{A}_{3}\mathbb{Z}+\mathbb{A}_{4}}, \,\mathbb{Z}\neq\frac{\mathbb{A}_{4}}{\mathbb{A}_{3}},$$
где $\mathbb{A}_{1}, \mathbb{A}_{2}, \mathbb{A}_{3}, \mathbb{A}_{4}$
-- фиксированние комплексные числа, а
$\mathbb{Z}=\{z_{k}\}_{k=1}^{\infty}$ -- комплексное переменное.

\textbf{2)} $W=\mathbb{Z}^{n}=\{z_{k}^{n}\}_{k=1}^{\infty}$, где $n$
-- натуральное число, -- cтепенная функция, голоморфна во всей
плоскости $\overline{\mathbb{C}}^{\infty}.$

\textbf{3)}
$W=\frac{1}{2}\left(\mathbb{Z}+\frac{1}{\mathbb{Z}}\right)$ --
функция Жуковского, голоморфна в
$\overline{\mathbb{C}}^{\infty}\setminus\Theta.$\vskip 1mm

\textbf{4)}
$\mathbb{P}_{n}(\mathbb{Z})=\sum\limits_{k=0}^{n}\mathbb{A}_{k}\mathbb{Z}^{k}$
-- полином, $\mathbb{Z}\in \mathbb{C}^{\infty}$. \vskip 1mm

\textbf{5)} $\frac{1}{\mathbb{Z}-\mathbb{Z}_{0}}$,
$\mathbb{Z}-\mathbb{Z}_{0}\in
\mathbb{C}^{\infty}\setminus\Theta.$\vskip 1mm

\textbf{6)}
$\exp\mathbb{Z}=\{e^{z_{k}}\}_{k=1}^{\infty}=\textbf{1}+\mathbb{Z}+\frac{1}{2}\mathbb{Z}^{2}+\ldots+\frac{1}{k!}\mathbb{Z}^{k}+\ldots
$, $\mathbb{Z}\in \mathbb{C}^{\infty}$.

\textbf{7)}
$(1-\mathbb{Z})^{\frac{1}{2}}=\textbf{1}-\frac{1}{2}\mathbb{Z}+\frac{1}{8}\mathbb{Z}^{2}-\ldots+\frac{\frac{1}{2}
\left(\frac{1}{2}-1\right)\ldots(\frac{1}{2}-k+1)}
{k!}\mathbb{Z}^{k}-\ldots,$ $\quad \mathbb{Z}\in \mathbb{
U}^{\infty}=\{\mathbb{Z}: \|z\|<1\}.$

{\bf 10. Полицилиндрическая теорема Римана об отображении в
$\overline{\mathbb{C}}^{\infty}$.}

Означим области гиперболичного типа как в пункте 10 \cite{32}.
Область $B\subset\overline{\mathbb{C}}$ называется областью
гиперболического типа, если $\partial B$ (граница B) -- связное
множество, содержащее более одной точки.

Пусть $0<\delta\leq1$ и
$\mathbb{A}=\{a_{k}\}_{k=1}^{\infty}\in\overline{\mathbb{C}}^{\infty}$.
Тогда
$\mathbb{B}=\mathbb{B}(\delta)=\mathbb{B}_{\delta}(\mathbb{A})=B_{1}\times
B_{2}\times\ldots\times
B_{n}\times\ldots\subset\overline{\mathbb{C}}^{\infty}$,
$\mathbb{A}\in\mathbb{B}_{\delta}(\mathbb{A})$, где каждая область
$B_{k}$ является областью гиперболического типа,
$\delta<r(B_{k},a_{k})<\frac{1}{\delta}$, $k=\overline{1,\infty}$.
При любом $0<\delta\leq1$ область
$\mathbb{B}(\delta)=\mathbb{B}_{\delta}(\mathbb{A})$ называется
конечной относительно $\mathbb{A}$ полицилиндрической областью
гиперболического типа.

Аналогично \cite{32} сформулируем утверждение непосредственно
вытекающее из классической теоремы Римана об отображении односвязной
области гиперболического типа на единичный круг.

\textbf{Теорема Римана.} \emph{Пусть
$\mathbb{A}\in\overline{\mathbb{C}}^{\infty}$ и $0<\delta\leq1$.
Тогда любая конечная относительно $\mathbb{A}$ полицилиндрическая
область
$\mathbb{B}=\mathbb{B}_{\delta}(\mathbb{A})\subset\overline{\mathbb{C}}^{\infty}$
гиперболического типа биголоморфно эквивалентна единичному поликругу
$\mathbb{U}^{\infty}=\{\mathbb{W}\in
\mathbb{C}^{\infty}:\|\mathbb{W}\|<1\}$.}

Пусть $\mathbb{B}=\mathbb{B}(\delta)=B_{1}\times
B_{2}\times\ldots\times B_{n}\times\ldots$ -- область указанная в
теореме Римана, $\mathbb{A}=\{a_{k}\}_{k=1}^{\infty}\in \mathbb{B}$,
$a_{k}\in B_{k}$, $k=\overline{1,\infty}$ и $w_{k}=f_{k}(z_{k})$ --
голоморфная в $B_{k}$ функция, однолистно и конформно отображающая
область $B_{k}$, $k=\overline{1,\infty}$ на единичный круг
$|w_{k}|<1$ так, что $f(a_{k})=0$, $f'(a_{k})>0$.

Тогда биголоморфное отображение
$\mathbb{F}_{\mathbb{B}}(\mathbb{Z})=\{f_{k}(z_{k})\}_{k=1}^{\infty},\quad
\mathbb{F}'_{\mathbb{B}}(\mathbb{Z})=\{f'_{k}\}_{k=1}^{\infty},$
удовлетворяет условиям нормировки
$$\mathbb{F}_{\mathbb{B}}(\mathbb{A})=\mathbb{O},\quad  \mathbb{F}_{\mathbb{B}}'(\mathbb{A})=\{f_{k}'(a_{k})\}_{k=1}^{\infty}
>\mathbb{O}$$
и будет единственним таким отображением на единичный поликруг. Тогда
обратное отображение к отображению
$\mathbb{F}_{\mathbb{B}}(\mathbb{A})$ является частично--конформным
отображением единичного поликруга.

Итак, в алгебре ${\mathbb{C}}^{\infty}$ норма определена равенством
$\|\mathbb{Z}\|:=|\mathbb{Z}|.$ Метрика (векторная) в
${\mathbb{C}}^{\infty}$ задается обычным образом:
$\rho(\mathbb{Z}_{1}, \mathbb{Z}_{2})=\|\mathbb{Z}_{1} -
\mathbb{Z}_{2}\|$. Назовем так определенные (векторные) норму и
метрику полицилиндрическими. Сходимость по полицилиндрической норме
равномерно по номерам задается соотношением
$\mathbb{Z}_{p}\underset{p\rightarrow\infty}\longrightarrow
0\Longleftrightarrow
\|\mathbb{Z}_{p}\|\underset{p\rightarrow\infty}\longrightarrow\mathbb{O}=(0,
0, \ldots 0, \ldots)\Longleftrightarrow
|z_{p}^{(k)}|\underset{p\rightarrow\infty}\rightrightarrows 0,
\forall k=\overline{1, \infty}.$ (Знак "$\rightrightarrows$"
обозначает равномерную сходимость по $k=\overline{1, \infty}$.)

\textbf{11. Приложения.} В связи с бесконечномерной теоремой Римана
об отображении рассмотрим полицилиндрический аналог известного
класса $S$ из теории однолистных функций \cite{b6} -- \cite{10}.

\textbf{\emph{Определение 13.}} Классом $\mathbb{S}^{(\infty)}$
назовем совокупность всех биголоморфных отображений единичного
поликруга $\mathbb{U}^{\infty}=\{\mathbb{Z}\in\mathbb{C}^{\infty}:
\|\mathbb{Z}\|<1\}$ вида
$\mathbb{F}(\mathbb{Z})=\{f_{k}(z_{k})\}_{k=1}^{\infty}$, где
$f_{k}\in S,$ $k=\overline{1, \infty}$,
$\mathbb{Z}=\{z_{k}\}_{k=1}^{\infty}\in \mathbb{U}^{\infty}$.

Ясно, что для $\mathbb{Z}\in \overline{\mathbb{U}}^{\infty}(r): =
\{\|\mathbb{Z}\|\leq r < 1\},$ $0<r<1,$ равномерно и абсолютно
сходится ряд

$$\mathbb{F}(\mathbb{Z})=\sum \limits_{k=1}^\infty
\mathbb{A}_{k}\mathbb{Z}^{k}=\sum\begin{pmatrix}
a_{k}^{(1)}\\
a_{k}^{(2)}\\
\vdots\\
a_{k}^{(n)}\\
\vdots\\
\end{pmatrix}\begin{pmatrix}
z_{1}\\
z_{2}\\
\vdots\\
z_{n}\\
\vdots\\
\end{pmatrix}^{k}=\begin{pmatrix}
\sum a_{k}^{(1)}z_{1}^{k}\\
\sum a_{k}^{(2)}z_{2}^{k}\\
\vdots\\
\sum a_{k}^{(n)}z_{n}^{k}\\
\vdots\\
\end{pmatrix}=\begin{pmatrix}f_{1}(z_{1})\\
f_{2}(z_{2})\\
\vdots\\
f_{n}(z_{n})\\
\vdots\\
\end{pmatrix}.$$

\textbf{\theor} Для произвольного отображения $\mathbb{F}\in
\mathbb{S}^{(\infty)}$ справедливо неравенство
$$\frac{\|\mathbb{Z}\|}{(1+\|\mathbb{Z}\|)^{2}}\leq \|\mathbb{F}(\mathbb{Z})\|\leq\frac{\|\mathbb{Z}\|}{(1-\|\mathbb{Z}\|)^{2}}\:,$$
где $\|\mathbb{Z}\|=r,$\quad $0\leq r<1$. \vskip 1mm

\textbf{\theor} Для произвольного отображения $\mathbb{F}\in
\mathbb{S}^{(\infty)}$ справедливо неравенство
$$\frac{\|1 - \mathbb{Z}\|}{(1+\|\mathbb{Z}\|)^{3}}\leq
\|\mathbb{F}'(\mathbb{Z})\|\leq\frac{\|1+\mathbb{Z}\|}{(1-\|\mathbb{Z}\|)^{3}}\:,$$
где $\|\mathbb{Z}\|=r,$\quad $0\leq r<1,$\quad $k=\overline{1,
\infty}$.

\textbf{\theor} (Теорема Бибербаха) Если $\mathbb{F}\in
\mathbb{S}^{(\infty)}$, тогда
$$|\mathbb{A}_{n}|\leq n\cdot \textbf{1} = n,$$
где $\mathbb{F}=\sum \limits_{k=1}^\infty
\mathbb{A}_{k}\mathbb{Z}^{k}$, $\textbf{1}=(1,1,\ldots, 1, \ldots)$.
Знак равенства в этом неравенстве достигается тогда и только тогда,
когда $\mathbb{F}=\{f_{k}\}_{k=1}^{\infty}$,
$f_{k}^{0}=z_{k}(1-e^{i\theta}z_{k})^{-2}$, $\theta_{k}\in[0,2\pi]$,
$k=\overline{1, \infty}$.

\textbf{12. Некоторые задачи о неналегающих областях в
$\overline{\mathbb{C}}^{\infty}$.}

Перенесем результат М.А. Лаврентьева о максимуме произведения
конформных радиусов двух неналегающих областей в пространство
$\overline{\mathbb{C}}^{\infty}$. Пусть $\mathbb{B}_{1}$ и
$\mathbb{B}_{2}$ -- две полицилиндрические области в
$\overline{\mathbb{C}}^{\infty}$, причем
$\mathbb{B}_{1}=B_{1}^{(1)}\times B_{2}^{(1)}\times
B_{3}^{(1)}\times...\times B_{n}^{(1)}\times\ldots$,
$\mathbb{B}_{2}=B_{1}^{(2)}\times B_{2}^{(2)}\times
B_{3}^{(2)}\times...\times B_{n}^{(2)}\times\ldots$, где
$B_{k}^{(s)}\subset \overline{\mathbb{C}}$, $k=\overline{1,
\infty}$, $s=\overline{1, 2}$ -- произвольные многосвязные области.
Пусть $\mathbb{A}_{1}\in\mathbb{B}_{1}$,
$\mathbb{A}_{2}\in\mathbb{B}_{2}$.

Обобщенным внутренним радиусом полицилиндрической области
$\mathbb{B}$ в точце $\mathbb{A}$ ($\mathbb{A}\in\mathbb{B}$) будем
называть неупорядоченный набор чисел
$\{r(B_{k},a_{k})\}_{k=1}^{\infty}$. Для удобства числа
$r(B_{k},a_{k})$ будем располагать в порядке координатных областей

$$\mathbb{R}(\mathbb{B},
\mathbb{A})=\begin{pmatrix}
r(B_{1},a_{1})\\
\ldots\\
r(B_{n},a_{n})\\
\ldots
\end{pmatrix}.$$

Рассмотрим класс пар полицилиндрических областей
$\Lambda(\mathbb{A}_{1},\mathbb{A}_{2})$, который содержит пары
полицилиндрических областей, для которых набор координатных областей
$\{B_{k}^{s}\}$, $k=\overline{1, \infty}$ при каждом фиксированном
$s=\overline{1, 2}$ является системой неналегающих областей на
$\overline{\mathbb{C}}$. На классе
$\Lambda(\mathbb{A}_{1},\mathbb{A}_{2})$ запишем следующий
функционал

\begin{equation}\label{1}
\mathbb{J}=\mathbb{R}(\mathbb{B}_{1},\mathbb{A}_{1})\cdot\mathbb{R}(\mathbb{B}_{2},\mathbb{A}_{2}).
\end{equation}

\textbf{Задача}. На классе $\Lambda(\mathbb{A}_{1},\mathbb{A}_{2})$
для фиксированных точек $\mathbb{A}_{1}=\mathbb{O}$ и
$\mathbb{A}_{2}=\infty$ определить максимум функционала (\ref{1}) и
описать все экстремали.

\textbf{\theor} (\emph{А.К. Бахтин, И.В. Денега}) На классе областей
$\Lambda(\mathbb{O},\infty)$ справедливо неравенство

\begin{center}
$\mathbb{R}(\mathbb{B}_{1},\mathbb{O})\cdot\mathbb{R}(\mathbb{B}_{2},\infty)\leq
\textbf{1},$
\end{center}
причем знак равенства в этом неравенстве достигается для областей:
$\mathbb{B}_{1}=\mathbb{U}^{\infty}$ и $\mathbb{B}_{2}$ -- образ
$\mathbb{B}_{1}$ при частично-- конформном отображении области
$\mathbb{B}_{1}$ с помощью функции
$\mathbb{F}=\frac{\textbf{1}}{\mathbb{Z}}$. $(\mathbb{O}=(0, 0,
...,0,...), \infty=(\infty, \infty,...,\infty,...))$.


\newpage

\begin{thebibliography}{99}
\bibitem{32} Бахтин А.К. {\sl Обобщение некоторых результатов теории однолистных функций на многомерные комплексные пространства} // Доп. НАН України. --
2011. -- №~3. -- С.~7~--~11.
\bibitem{b1} Шабат Б. В. {\sl Введение в комплексный
анализ}, Ч. І. -- М.:«Наука», 1976. -- 320 с.
\bibitem{b2} Шабат Б. В. {\sl Введение в комплексный анализ}, Ч. ІІ. -- М.:«Наука», 1976. -- 400 с.
\bibitem{b3}  Фукс Б. В.  {\sl
Введение в теорию аналитических функций многих комплексных
переменных}, Физматгиз, 1962. -- 420 с.
\bibitem{b4} Фукс Б. В. {\sl Специальные главы теории аналитических функций многих комплексных
переменных}, Физматгиз, 1963. -- 428 с.
\bibitem{b5} Чирка Е. М. {\sl Комплексные аналитические множества.} -- М.:«Наука», 1985. -- 272 с.
\bibitem{b6}  Голузин Г. М. {\sl Геометрическая теория функций комплексного
переменного}. -- М: «Наука», 1966. -- 628 с.
\bibitem{heim} Хейман В К. {\sl Многолистные функции. - М.: Изд-во иностр.
лит}., 1960. -- 180 с.
\bibitem{3} Дубинин В. Н. {\sl Метод симметризации в задачах о неналегающих
областях} // Мат. сб. -- 1985. -- {\bf 128}, №~1. -- С.~110~--~123.
\bibitem{10}  Дубинин В. Н.{\sl Емкости конденсаторов и симметризация в
геометрической теории функций комплексного переменного.}
// Владивосток "Дальнаука" ДВО РАН -- 2009. -- 390с.
\bibitem{b7} Кантор И. Л., Солодовников А. С. {\sl Гиперкомплексные числа}. -- М.:«Наука», 1973. -- 143 с.
\bibitem{b8} Б. Л. ван дер ВАРДЕН. {\sl Алгебра} -- М.:«Наука», 1976. -- 648 с.
\bibitem{b10} Рудин У. {\sl Функциональный анализ}, -- Изд. "Мир", Москва. -- 1975. -- 449 с.

\bibitem{b9} Шилов Г. Е. {\sl Математический анализ. Конечномерные линейные пространства}, -- М.:«Наука». -- 1969. -- 432 с.

\end{thebibliography}
\end{document}